\documentclass[11pt,a4paper]{article}

\usepackage{graphics}
\usepackage{amsmath,amsthm,amssymb}

\makeatletter
\thm@notefont{\fontseries\mddefault\upshape}
\thm@headfont{\scshape}
\makeatother

\theoremstyle{plain}
\newtheorem{theorem}{Theorem}
\newtheorem*{theorem*}{Theorem}

\newtheorem*{proposition*}{Proposition}
\newtheorem{lemma}[theorem]{Lemma}
\newtheorem*{lemma*}{Lemma}

\newtheorem*{corollary*}{Corollary}
\theoremstyle{remark}

\newtheorem*{remark*}{Remark}
\theoremstyle{definition}
\newtheorem*{example*}{Example}
\newtheorem*{acknowledgment}{Acknowledgment}

\def\SZ{\mathbf{Z}}

\def\SC{\mathbf{C}}
\def\SP{\mathbf{P}}
\def\Aut{\mathop{\rm Aut}\nolimits}

\title{%
The projective characterization of elliptic plane curves which have one place at infinity}
\author{%
Keita Tono}
\date{}

\begin{document}
\maketitle
\def\thefootnote{}
\def\thefootnote{1}
\begin{abstract}
In this paper
we consider smooth affine elliptic plane curves
having one place at infinity.
We identify them with elliptic projective plane curves
having only one cusp as their singular points
and meeting with the line at infinity only at the cusp.
We characterize such curves by
the self-intersection number of the strict transform of them
via the minimal embedded resolution of their cusp.
Furthermore,
we prove that the self-intersection number of them
is the maximum value among those of
all the elliptic plane curves having only one cusp.
\end{abstract}
\section{Introduction}
Let $C$ be a curve on $\SP^2=\SP^2(\SC)$.
A singular point of $C$ is said to be a \emph{cusp}
if it is a locally irreducible singular point.
We say that $C$ is \emph{cuspidal} (resp.~\emph{unicuspidal})
if $C$ has only cusps (resp.~one cusp) as its singular points.
Suppose $C$ is unicuspidal.
We denote by $C'$ the strict transform of $C$
via the minimal embedded resolution of the cusp of $C$.
In view of \cite{am,su},
we say that $C$ is of
\emph{Abhyankar-Moh-Suzuki type} (\emph{AMS type}, for short)
if there exists a line $L$ such that $C\cap L=\{\text{the cusp}\}$.
By regarding $L$ as the line at infinity,
we identify such curves with
smooth affine plane curves having one place at infinity.
In \cite{yoshihara2},
it was proved that
a rational unicuspidal plane curve $C$ is of AMS type
if and only if $(C')^2\ge2$.
In a similar manner to that as in the rational case,
we can characterize
elliptic unicuspidal plane curves $C$ of AMS type by $(C')^2$.
Namely, the purpose of this paper is to prove the following:
\begin{theorem}\label{thm:ams}
If $C$ is an elliptic unicuspidal plane curve,
then $(C')^2\le6$.
The equality holds if and only if $C$ is of AMS type.
\end{theorem}

We next determine the maximum value of $(C')^2$ for
elliptic unicuspidal plane curves
of non-AMS type.
\begin{theorem}\label{thm:nonams}
If $C$ is an elliptic unicuspidal plane curve of non-AMS type,
then $(C')^2\le3$.
The equality holds if and only if
there exist an irreducible conic $C_2\subset\SP^2$ and
a birational map $f:\SP^2\rightarrow\SP^2$ such that
$C\cap C_2=\{\text{the cusp}\}$,
$f|_{\SP^2\setminus C_2}\in\Aut(\SP^2\setminus C_2)$
and $f(C)$ is a smooth cubic curve.
Furthermore, we can show the existence of such curves $C$ with $(C')^2=3$.
\end{theorem}

In Section~\ref{sec:pf},
we prove Theorem~\ref{thm:ams} except the ``if'' part.
After resolving the cusp of an elliptic unicuspidal plane curve $C$,
we perform additional blowings-up to make $\mathrm{\Phi}_{|C'|}$
an elliptic fibration.
Then we analyze its structure
and determine  possible types of its singular fibers.
We show that if $(C')^2=6$, then the fibration has a singular fiber
of type $\mathrm{II}^{\ast}$.
In Section~\ref{sec:pfif},
we prove the ``if'' part of Theorem~\ref{thm:ams}
by using the following property of
elliptic unicuspidal plane curves of AMS type.
\begin{theorem}[{\cite[Theorem 8.7]{ao}, cf.~\cite{miyanishi}}]\label{thm:ao}
Let $C$ be an elliptic unicuspidal plane curve of AMS type
and
$L$ a line such that $C\cap L=\{\text{the cusp}\}$.
Then there exists a birational map $f:\SP^2\rightarrow\SP^2$ such that
$f|_{\SP^2\setminus L}\in\Aut(\SP^2\setminus L)$
and $f(C)$ is a smooth cubic curve.
\end{theorem}

In Section~\ref{sec:nonams},
we prove Theorem~\ref{thm:nonams}.
Similar to Section~\ref{sec:pf},
we consider the elliptic fibration associated with $|C'|$.
We prove that if $(C')^2=3$, then the fibration has a singular fiber
of type $\mathrm{I}_4^{\ast}$.
In order to show the existence of curves $C$ with $(C')^2=3$,
we give the defining equation of the cubic pencil
associated with the elliptic fibration.
%
\section{%
Proof of Theorem~\ref{thm:ams}}\label{sec:pf}
Let $C$ be an elliptic unicuspidal plane curve
and $P$ the cusp of $C$.
Let $\sigma:V\rightarrow\SP^2$ denote the minimal embedded resolution
of $P$.
That is,
$\sigma$ is the composite of the shortest sequence
of blowings-up such that
the strict transform $C'$ of $C$ intersects $\sigma^{-1}(P)$ transversally.
The dual graph of $D:=\sigma^{-1}(C)$ has the following shape,
where $g\ge1$ and all $A_i$, $B_i$ are not empty.
\begin{center}
\includegraphics{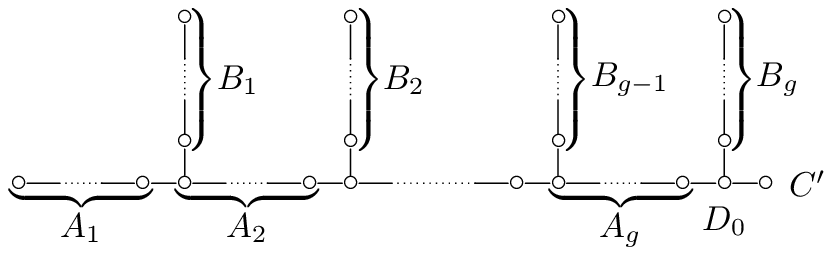}
\end{center}
Here $D_0$ is the exceptional curve of the last blowing-up
and
$A_1$ contains that of the first one.
The morphism $\sigma$ contracts $A_g+D_0+B_g$ to a ($-1$)-curve $E$,
$A_{g-1}+E+B_{g-1}$ to a ($-1$)-curve and so on.
Every irreducible component $E$ of
$A_i$ and $B_i$ is a smooth rational curve with $E^2<-1$.
Each $A_i$ contains an irreducible component $E$
such that $E^2<-2$.
Cf.~\cite{bk,masa}.
We give weights to $A_1,\ldots,A_g,B_1,\ldots,B_g$
in the usual way.
We also give the graphs $A_1,\ldots,A_g$ (resp.~$B_1,\ldots,B_g$)
the direction from the left-hand side to the right
(resp.~from the bottom to the top) in the above figure.
Let $D_1$ (resp.~$D_2$) denote
the first irreducible component of $B_g$ (resp.~the last one of $A_g$)
with respect to the direction.
%

Suppose that $n:=(C')^2\ge3$.
Perform ($n-1)$-times of blowings-up $\tau_0:W_0\rightarrow V$
over $C'\cap D_0$ in the following way,
where
$\ast$ (resp.~$\bullet$) denotes
a ($-1$)-curve (resp.~($-2$)-curve)
and $E_i$ is the exceptional curve of the $i$-th blowing-up.
We use the same notation $C'$.
\begin{center}
\scalebox{1}{\includegraphics{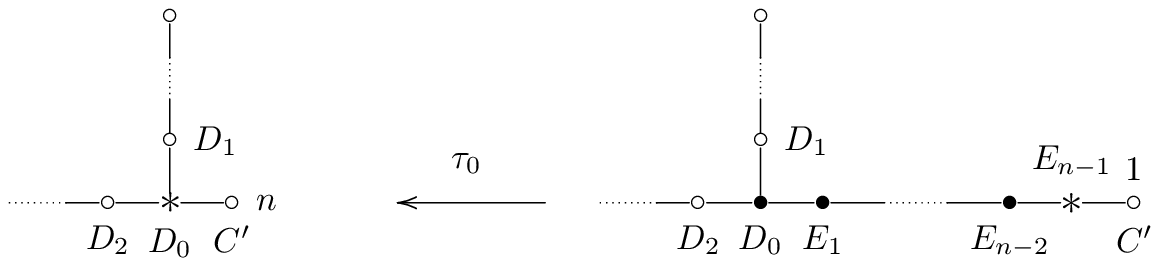}}
\end{center}

On $W_0$, there exists an exact sequence:
\[
0\longrightarrow H^0(\mathcal{O}_{W_0})
\longrightarrow H^0(\mathcal{O}_{W_0}(C'))
\longrightarrow H^0(\mathcal{O}_{C'}(C'))
\longrightarrow H^1(\mathcal{O}_{W_0})=0.
\]
Because $(C')^2=1$ on $W_0$,
we have $h^0(\mathcal{O}_{C'}(C'))=1$.
This means that $|C'|$ is a pencil on $W_0$ having a single base point $Q$.
We have either $\{Q\}=C'\cap E_{n-1}$ or $Q\in C'\setminus E_{n-1}$.
Let $\tau_1:W\rightarrow W_0$ denote the blowing-up  at $Q$
and $E_n$ its exceptional curve.
Put $\tau=\tau_0\circ\tau_1$.
The pencil $|C'|$ does not have base points on $W$.
The morphism $\mathrm{\Phi}_{|C'|}:W\rightarrow\SP^1$
is an elliptic fibration.
The curve $C'$ is a nonsingular fiber.
The curve $E_{n}$ is a 1-section of $\mathrm{\Phi}_{|C'|}$.
Namely, we have $E_{n}C'=1$.
If $Q\in C'\setminus E_{n-1}$,
then $E_{n-1}$ is also a 1-section.
Put $F_0:=\tau^{-1}(D)-E_{n}-C'$ if $\{Q\}=C'\cap E_{n-1}$,
$F_0:=\tau^{-1}(D)-E_n-E_{n-1}-C'$ otherwise.
The divisor $F_0$
is contained in a fiber $\overline{F}_0$ of $\mathrm{\Phi}_{|C'|}$.
\begin{lemma}\label{lem1}
Let $F$ be a fiber of $\mathrm{\Phi}_{|C'|}$ other than $\overline{F}_0$.
If $\{Q\}=C'\cap E_{n-1}$,
then $F$ is irreducible.
Otherwise $F$ contains at most two irreducible components.
\end{lemma}
\begin{proof}
If $\{Q\}=C'\cap E_{n-1}$,
then $W\setminus(F_0\cup E_n\cup C')=\SP^2\setminus C$.
The latter does not contain a complete curve.
This implies that each irreducible component of $F$
intersects $E_n$.
Since $E_n$ is a 1-section,
$F$ must be irreducible.
Similarly,
if $Q\in C'\setminus E_{n-1}$,
each irreducible component of $F$
intersects $E_{n-1}$ or $E_n$.
It follows that $F$ has at most two irreducible components.
\end{proof}
Let $\varphi:W\rightarrow X$
be successive contractions of ($-1$)-curves
in the singular fibers of $\mathrm{\Phi}_{|C'|}$
such that the fibration
$p:=\mathrm{\Phi}_{|C'|}\circ\varphi^{-1}:X\rightarrow\SP^1$
is relatively minimal.
We will use the following facts about
relatively minimal elliptic rational surfaces
(\cite[Lemma~2.7, 10.2, Theorem~10.3]{shioda1}).
\begin{lemma}\label{lem2}
The following assertions hold on $X$.
\begin{enumerate}
\item[(i)]
Every 1-section is a ($-1$)-curve.
\item[(ii)]
The Picard number $\rho(X)$ of $X$ is equal to $10$.
\item[(iii)]
$\displaystyle \sum_{F}(r(F)-1)\le8$,
where $F$ runs over all fibers of $p$
and $r(F)$ denotes the number of irreducible components of $F$.
\end{enumerate}
\end{lemma}

We prove the following three lemmas
by using the above Lemmas together with
the list of all possible singular fibers
of relatively minimal elliptic fibrations given by
Kodaira (\cite[Theorem 6.2]{kodaira}).
\begin{lemma}\label{lem3}
The following assertions hold.
\begin{enumerate}
\item[(i)]
$\varphi$ does not contract any irreducible curve
which is not contained in $\overline{F}_0$.
\item[(ii)]
If $\{Q\}=C'\cap E_{n-1}$ (resp.~$Q\in C'\setminus E_{n-1}$),
then
$\varphi(D_0),\varphi(E_1),\ldots$, $\varphi(E_{n-1})$
(resp.~$\varphi(D_0),\varphi(E_1),\ldots,\varphi(E_{n-2})$)
are ($-2$)-curves.
\item[(iii)]
$\varphi$ does not contract $D_1,D_2$.
\end{enumerate}
\end{lemma}
\begin{proof}
The assertion (i) follows from Lemma~\ref{lem1} and Lemma~\ref{lem2} (i).
The assertion (iii) follows from (ii).
We show (ii).
We only prove the assertion for the case in which $\{Q\}=C'\cap E_{n-1}$.
The morphism $\varphi$ does not contract $E_{n-1}$.
If $\varphi$ contracts a curve meeting with $E_{n-1}$
then $\varphi(E_{n-1})^2>-2$.
It follows that $\varphi(E_{n-1})$ coincides with $\varphi(\overline{F}_0)$
and has a double point.
We have $0=\varphi(E_{n-1})^2\ge E_{n-1}^2+4=2$, which is a contradiction.
Thus $\varphi$ does not contract $E_{n-2}$.
So, the fiber $\varphi(\overline{F}_0)$ is reducible.
We infer that
$\varphi(\overline{F}_0)$ consists of ($-2$)-curves.
It turns out that $\varphi$ does not contract $D_0,E_1,\ldots,E_{n-3}$.
\end{proof}
\begin{lemma}\label{lem4}
The following assertions hold.
\begin{enumerate}
\item[(i)]
The morphism $\varphi$ is not the identity.
\item[(ii)]
If $\{Q\}=C'\cap E_{n-1}$
(resp.~$Q\in C'\setminus E_{n-1}$),
then
$r(\overline{F}_0)-r(F_0)=1$
(resp.~$r(\overline{F}_0)-r(F_0)=2$).
\item[(iii)]
$r(\varphi(\overline{F}_0))=9$.
\end{enumerate}
\end{lemma}
\begin{proof}
(i)
The morphism $\varphi$ is not the identity
because $F_0$ contains an irreducible component $E$ with $E^2\le-3$.

(ii)
Since $F_0$ does not contain a ($-1$)-curve,
we have $1\le r(\overline{F}_0)-r(F_0)$.
The number of the blowings-down of $\varphi$ is equal to
$\rho(W)-\rho(X)=r(D)+n-10$.
By Lemma~\ref{lem2} (iii),
we get $9\ge r(\varphi(\overline{F}_0))=r(\overline{F}_0)-(r(D)+n-10)$.
Thus $r(\overline{F}_0)\le r(D)+n-1$.
Since $r(F_0)=r(D)+n-2$ (resp.~$r(F_0)=r(D)+n-3$),
we have $1\le r(\overline{F}_0)-r(F_0)\le1$
(resp.~$1\le r(\overline{F}_0)-r(F_0)\le2$).
Suppose that $Q\in C'\setminus E_{n-1}$
and $r(\overline{F}_0)-r(F_0)=1$.
Let $E_0$ denote the ($-1$)-curve
which is contracted by the first blowing-down of $\varphi$.
Because $F_0\cap E_{n}=\emptyset$,
$E_0$ must intersect $E_n$, which contradicts Lemma~\ref{lem2} (i).

(iii)
We have $r(\overline{F}_0)=r(D)+n-1$.
Thus
$r(\varphi(\overline{F}_0))=r(\overline{F}_0)-(r(D)+n-10)=9$.
\end{proof}

Let $E_0$ denote the ($-1$)-curve
which is contracted by the first blowing-down of $\varphi$.
If $\{Q\}=C'\cap E_{n-1}$ then $\overline{F}_0=F_0+E_0$.
If $Q\in C'\setminus E_{n-1}$,
write $\overline{F}_0$ as
$\overline{F}_0=F_0+E_0+E_0'$,
where $E_0'$ is an irreducible curve.
The curve $E_0'$ must intersect $E_n$,
because $F_0+E_0$ does not.
By Lemma~\ref{lem2} (i), $\varphi$ does not contract $E_0'$.
\begin{lemma}\label{lem5}
The following assertions hold.
\begin{enumerate}
\item[(i)]
$\varphi(\overline{F}_0)$ is of type $\mathrm{II}^{\ast}$ or $\mathrm{I}_4^{\ast}$.
We have $E_0F_0=E_0'F_0=1$.
\item[(ii)]
If $\{Q\}=C'\cap E_{n-1}$,
then $n=6$.
The fiber $\varphi(\overline{F}_0)$ is of type $\mathrm{II}^{\ast}$.
\item[(iii)]
If $Q\in C'\setminus E_{n-1}$,
then $n=3$.
The fiber $\varphi(\overline{F}_0)$ is of type $\mathrm{I}_4^{\ast}$.
\end{enumerate}
\end{lemma}
\begin{proof}
(i)
By Lemma~\ref{lem4} (iii),
the fiber $\varphi(\overline{F}_0)$ is of type
$\mathrm{II}^{\ast}$, $\mathrm{I}_4^{\ast}$ or $\mathrm{I}_9$.
By Lemma~\ref{lem3},
the curve
$\varphi(D_0)$ is a branching component of $\varphi(\overline{F}_0)$.
Namely,
$\varphi(D_0)(\varphi(\overline{F}_0)-\varphi(D_0))\ge3$.
Thus $\varphi(\overline{F}_0)$ is not of type $\mathrm{I}_9$.
It follows that $E_0F_0=E_0'F_0=1$.

(ii)
By Lemma~\ref{lem3} and Lemma~\ref{lem4} (iii),
we have $n\le7$.
The dual graph of $\varphi(F_0)$ contains the following graph.
\begin{center}
\scalebox{1}{\includegraphics{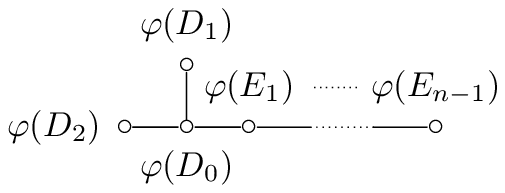}}
\end{center}
Suppose that $\varphi(E_i)$ is a branching component of
$\varphi(\overline{F}_0)$ for some $i<n$.
Let $E$ denote the irreducible component of $\varphi(\overline{F}_0-E_i-D_0)$
meeting with $\varphi(E_i)$.
Let $E'$ denote the strict transform of $E$ via $\varphi$.
Since $E'\ne E_0$, we have $E'\subset F_0$.
Because $E'$ does not intersect $E_i$,
$\varphi$ must contract a curve meeting with $E_i$,
which is a contradiction.
Hence $\varphi(E_1+\cdots+E_{n-1})$ does not contain a branching component.
Since $n\ge3$, $\varphi(\overline{F}_0)$ must be of type $\mathrm{II}^{\ast}$.
Because $\varphi(E_{n-1})$ intersects the section,
the coefficient of it in $p^{\ast}(p(\varphi(\overline{F}_0)))$
is equal to $1$.
This shows $n=6$.

(iii)
Since $\varphi$ does not contract $E_0'$,
we have $n\le7$ by Lemma~\ref{lem3} and Lemma~\ref{lem4} (iii).
If $\varphi(\overline{F}_0)$ is of type $\mathrm{II}^{\ast}$,
then $\varphi(E_n)$ must intersect $\varphi(E_{n-2})$,
which is impossible.
Thus $\varphi(\overline{F}_0)$ is of type $\mathrm{I}_4^{\ast}$.
The dual graph of $\varphi(\overline{F}_0-E_0')$ contains
the above graph with $\varphi(E_{n-1})$ being replaced with $\varphi(E_{n-2})$.
The divisor $\varphi(E_1+\cdots+E_{n-2})$ contains at most one
branching component of $\varphi(\overline{F}_0)$.
If it does,
then $n=7$ and $E_0'$ must intersect $E_{4}$.
Otherwise we have $n=3$.
Suppose $n=7$.
Let $\psi:X\rightarrow\SP^2$ denote the contraction of
$\varphi(E_7)+\varphi(E_6)+\varphi(\overline{F}_0)-\varphi(E_0')-\varphi(D_1)$.
Then we have $\psi(\varphi(D_1))^2=0$, which is absurd.
\end{proof}

Suppose $n=6$.
Let $C_0$ be an irreducible component of $\varphi(\overline{F}_0)$
whose position in the dual graph of $\varphi(\overline{F}_0)$
is illustrated in the following figure.
\begin{center}
\scalebox{1}{\includegraphics{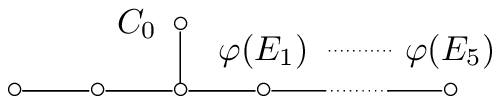}}
\end{center}
Let $\psi:X\rightarrow\SP^2$ denote the contraction of
$\varphi(E_6)+\varphi(\overline{F}_0)-C_0$.
We have $\psi(C_0)^2=1$.
Put $C_1=\sigma(\tau(E_0))$,
$h=\psi\circ\varphi\circ\tau^{-1}\circ\sigma^{-1}:\SP^2\rightarrow\SP^2$.
The morphism $h|_{\SP^2\setminus C_1}:\SP^2\setminus C_1\rightarrow\SP^2\setminus\psi(C_0)$ is an isomorphism.
Since
$\SZ/(\deg C_1)\SZ\cong H_1(\SP^2\setminus C_1,\SZ)\cong H_1(\SP^2\setminus\psi(C_0),\SZ)\cong\{0\}$,
we conclude that $C_1$ is a line (cf.~\cite[Proposition~4.1.3]{dimca}).
We completed the proof of Theorem~\ref{thm:ams} except the ``if'' part.
\section{%
Proof of Theorem~\ref{thm:ams} --- continued}\label{sec:pfif}
We use the following theorem to prove
the ``if'' part of Theorem~\ref{thm:ams}.
\begin{theorem}[{\cite[Lemma 4.4]{yoshihara}}]\label{thm:y}
Let $L\subset\SP^2$ be a line
and
$h:\SP^2\rightarrow\SP^2$ a birational map such that
$h\not\in\Aut\SP^2$ and
$h|_{\SP^2\setminus L}\in\Aut(\SP^2\setminus L)$.
Let $\sigma_h:V\rightarrow\SP^2$ denote the minimal resolution
of the base points of $h$.
Then the weighted dual graph of
$\sigma_h^{-1}(L)$ has the following shape,
where
$k\ge1$,
$n_i\ge0$,
$E_h$ is the exceptional curve of the last blowing-up.
\end{theorem}
\begin{center}
\scalebox{1}{\includegraphics{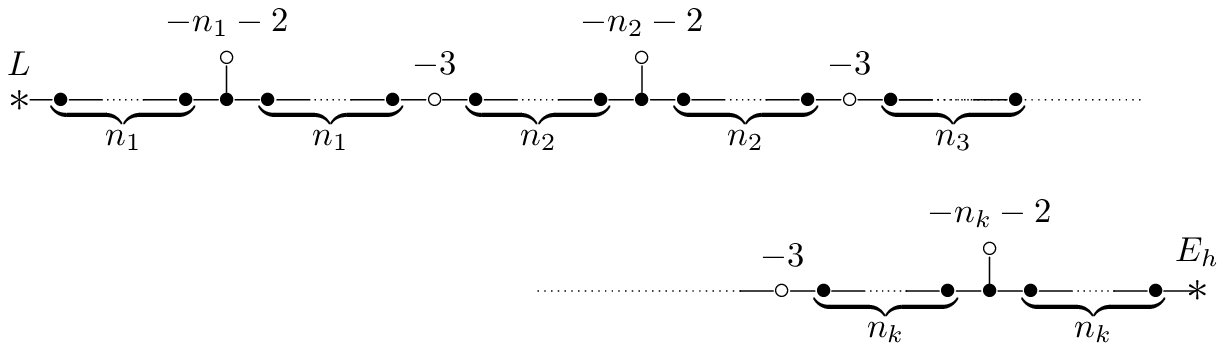}}
\end{center}

The morphism $h\circ\sigma_h$ contracts
the above graph in the following way,
where each marked subgraph is contracted to a point.
After the contraction,
$h\circ\sigma_h$ maps the image of $E_h$ under the contraction to $L$.
\begin{center}
\scalebox{1}{\includegraphics{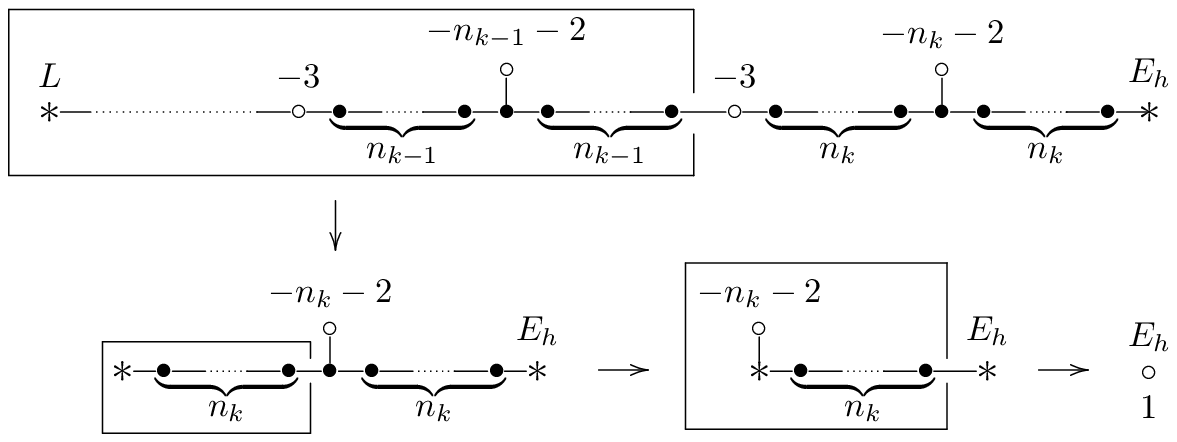}}
\end{center}

Let $C$ be an elliptic unicuspidal plane curve of AMS type.
Let $L$ and $f$ have the same meaning as in Theorem~\ref{thm:ao}.
We apply Theorem~\ref{thm:y} to $L$ and $f$.
On $\SP^2$, $f^{-1}$ has a single base point $Q\in L$.
Suppose that $Q\in L\setminus f(C)$.
Then
the center of every blowing-up of $f\circ\sigma_f:V\rightarrow\SP^2$
is not on $f(C)$.
This means that
the strict transform of $f(C)$ via $f\circ\sigma_f$
intersects $E_f$ in the same way as $f(C)$ does $L$.
Let $\sigma':V'\rightarrow V$ denote
the 3-times of blowings-up such that
the strict transform $C'$ of $f(C)$ on $V'$
intersects $(\sigma_f\circ\sigma')^{-1}(L)$
transversally.
We have $(C')^2=6$.
The center of every blowing-up of $\sigma_f\circ\sigma'$
is on $C$.
Thus $\sigma_f\circ\sigma'$ is
the minimal embedded resolution of the cusp of $C$.
The following figure illustrates
the weighted dual graph of $(\sigma_f\circ\sigma')^{-1}(C)$ near $C'$.
\begin{center}
\scalebox{1}{\includegraphics{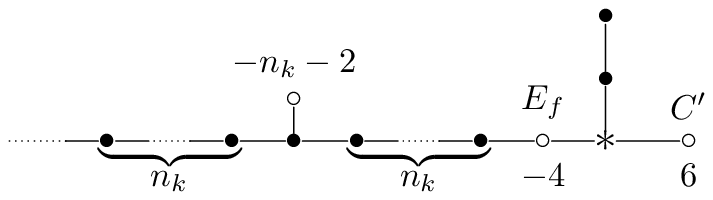}}
\end{center}

Suppose that $\{Q\}=L\cap f(C)$.
The center of every blowing-up of $f\circ\sigma_f:V\rightarrow\SP^2$
except the first two
is not on $f(C)$.
By arguments similar to those in the previous case,
we conclude that $(C')^2=6$,
where $C'$ is the strict transform of $C$ via
the minimal embedded resolution $\sigma:V\rightarrow\SP^2$ of the cusp.
The following figure illustrates
the weighted dual graph of $\sigma^{-1}(C)$ near $C'$.
\begin{center}
\scalebox{1}{\includegraphics{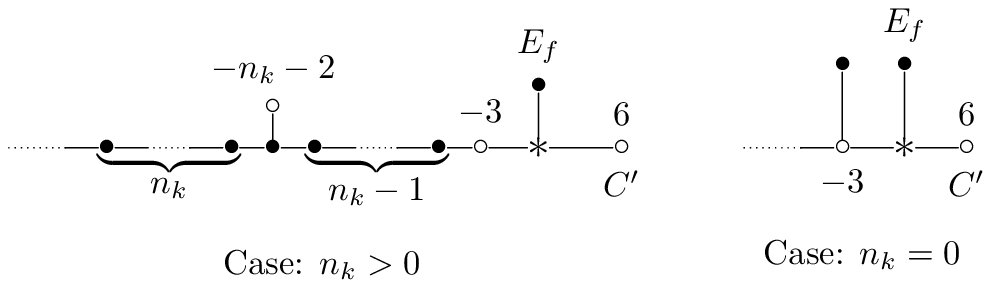}}
\end{center}
\section{%
Proof of Theorem~\ref{thm:nonams}}\label{sec:nonams}
Let the notation and the assumptions be as in Section~\ref{sec:pf},
where we assumed that $n\ge3$.
Assume in addition that $n\le5$.
By Lemma~\ref{lem5}, we must have $n=3$.
This proves the first inequality.
Moreover $\varphi(\overline{F}_0)$ is of type $\mathrm{I}_4^{\ast}$.
Since the coefficient of $\varphi(E_0')$
in $p^{\ast}(p(\varphi(\overline{F}_0)))$ is equal to one,
it follows that the weighted dual graph of
$\varphi(\overline{F}_0+E_2+E_3+E_0'+C')$
has the following shape.
\begin{center}
{\scalebox{1}{\includegraphics{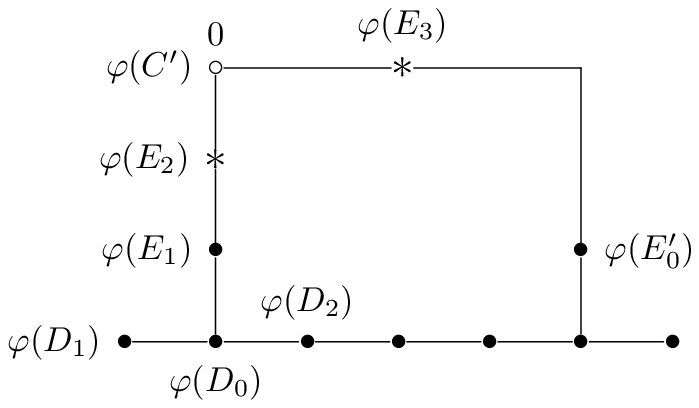}}}
\end{center}

We perform the blowings-down $\psi:X\rightarrow\SP^2$
as illustrated in the following figure,
where each marked subgraph is contracted to a point.
\begin{center}
{\scalebox{1}{\includegraphics{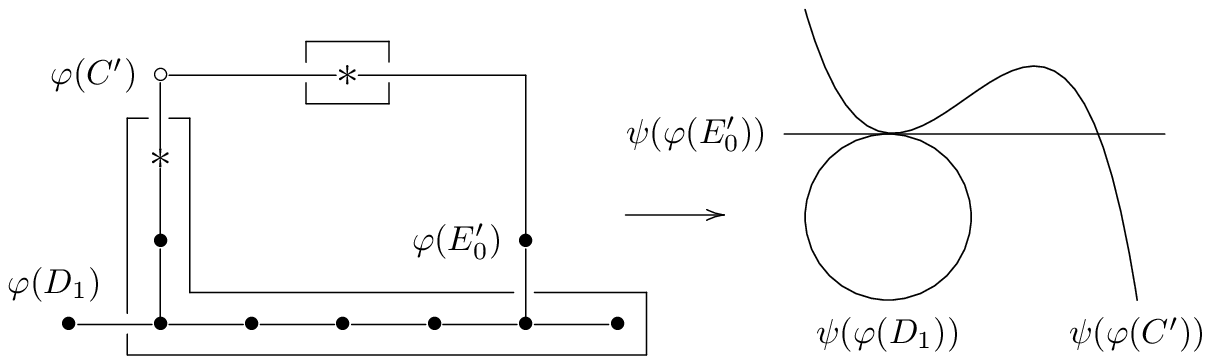}}}
\end{center}

The curve $L:=\psi(\varphi(E_0'))$ is a line.
The curves $C_2:=\psi(\varphi(D_1))$,
$C_3:=\psi(\varphi(C'))$ are a conic and a cubic curve, respectively.
We have
\renewcommand{\theequation}{$\ast$}
\begin{equation}\label{eqn1}
\text{$C_3\cdot L=R_1+2R_2$ and $C_3\cdot C_2=6R_2$}
\end{equation}
for some distinct points $R_1,R_2$.
Put $C_1=\sigma(\tau(E_0))$,
$f=\psi\circ\varphi\circ\tau^{-1}\circ\sigma^{-1}$.
The morphism $f|_{\SP^2\setminus C_1}:\SP^2\setminus C_1\rightarrow\SP^2\setminus C_2$ is an isomorphism.
This means that $C_1$ is an irreducible conic.
Hence the ``only if'' part is proved.
We prove the converse.
\begin{lemma}[{cf.~\cite[Lemma 4.4]{yoshihara}}]\label{lem:y2}
Let $C_2\subset\SP^2$ be an irreducible conic
and
$h:\SP^2\rightarrow\SP^2$ a birational map such that
$h\not\in\Aut\SP^2$ and
$h|_{\SP^2\setminus C_2}\in\Aut(\SP^2\setminus C_2)$.
Let $\sigma_h:V\rightarrow\SP^2$ denote the minimal resolution
of the base points of $h$.
Then the weighted dual graph of
$\sigma_h^{-1}(C_2)$ has the following shape,
where
$k\ge1$,
$n_i\ge0$,
$E_h$ is the exceptional curve of the last blowing-up.
Starting with the contraction of $C_2$,
$h\circ\sigma_h$ contracts the graph and maps the image of $E_h$
under the contraction to $C_2$.
\begin{center}
{\scalebox{1}{\includegraphics{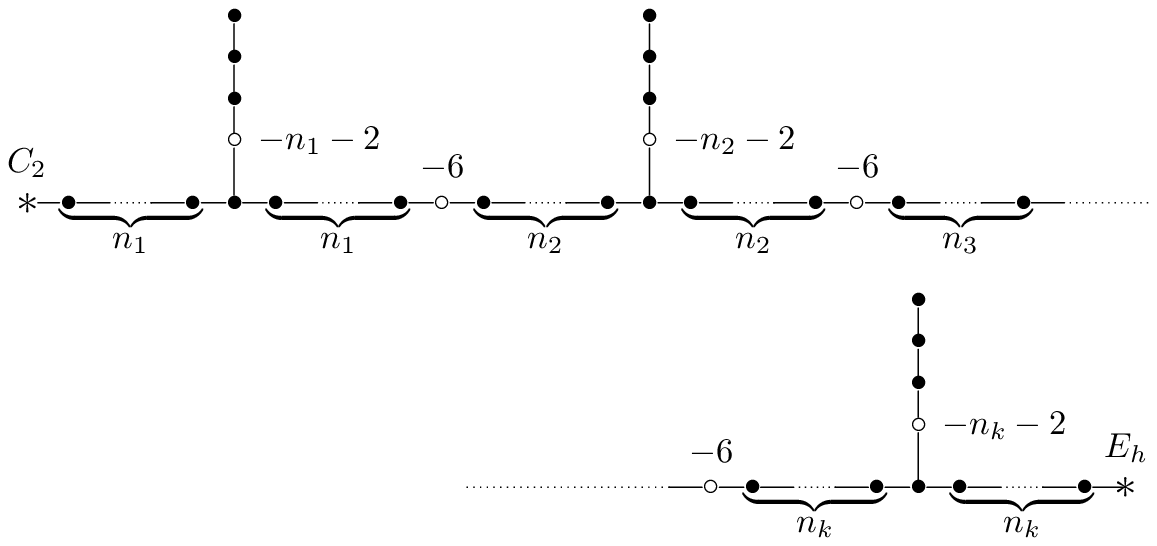}}}
\end{center}
\end{lemma}
\begin{proof}
We only give a sketch
because the proof is similar to that of {\cite[Lemma 4.4]{yoshihara}}.
Let $\sigma_i:V_i\rightarrow V_{i-1}$
denote the $i$-th blowing-up of $\sigma_h$,
$Q_{i-1}$ its center and $E_i$ its exceptional curve
($i=1,\ldots,r$, $V_0=\SP^2$).
We use the same symbols for the strict transforms
of $C_2$ and $E_i$.
We use the following facts ({\cite[Lemma 2.4, Remark 2.5]{yoshihara}}).
\begin{lemma}\label{lem:y3}
The following assertions hold for $i<r$.
\begin{enumerate}
\item[(i)]
We have $Q_i\in E_i$.
On $V_r$, we have $E_i^2\le-2$ and $C_2^2=-1$.
\item[(ii)]
If $C_2^2\ge0$ (resp.~$C_2^2=-1$) on $V_i$,
then $Q_i\in E_i\cap C_2$ (resp.~$Q_i\in E_i\setminus C_2$).
\item[(iii)]
$h\circ\sigma_h$ first contracts $C_2$,
then contracts $E_1,\ldots,E_{r-1}$, not necessarily in this order.
It does not contract $E_r$.
\end{enumerate}
\end{lemma}
We have $r\ge5$,
$Q_i\in E_i\cap C_2$ for $i<5$
and
$Q_i\in E_i\setminus C_2$ for $i\ge5$.
If $r=5$, then $h\circ\sigma_h$ cannot contract
the remaining curves after the contraction of $C_2$.
Thus $r\ge6$.
For the same reason, we infer $Q_5\in E_5\setminus E_4$
if $r=6$.
We have $k=1$, $n_1=0$ in this case.

Suppose $r\ge7$.
We have either
$\{Q_5\}=E_5\cap E_4$ or $Q_5\in E_5\setminus E_4$.
Suppose $\{Q_5\}=E_5\cap E_4$.
If $Q_6\in E_6\cap E_5$,
then $E_5^2\le-3$ on $V_r$.
This means that
none of $E_1,\ldots,E_{r-1}$ is a ($-1$)-curve
after the contraction of $C_2$.
Thus $Q_6\in E_6\setminus E_5$.
Let $\Gamma_i$ denote
the preimage of $C_2$ under $\sigma_1\circ\cdots\circ\sigma_i$.
Suppose that $Q_i$ is a node of $\Gamma_i$ for all $i\ge5$.
For the same reason as above, we have $Q_i\in E_i\cap E_4$ for all $i\ge5$.
After the contraction of $C_2,E_5,\ldots,E_{r-1}$,
$h\circ\sigma_h$ cannot contract the remaining curves.
Hence
there exists $n\ge1$ such that
$Q_{5+i}$ is a node of $\Gamma_{5+i}$ for $i=0,\ldots,n-1$
but $Q_{5+n}$ is not.
We set $n=0$ when $Q_5\in E_5\setminus E_4$.

Suppose that $Q_{6+i}$ is not a node of $\Gamma_{6+i}$ for all $i\ge n$.
We infer that $h\circ\sigma_h$ contracts
$C_2,E_5,\ldots,E_{r-1},E_4,E_3,E_2,E_1$
in this order.
After the contraction of $E_{r-1}$, $E_4$ must be a ($-1$)-curve.
Thus $r=2n+6$.
We have $k=1$, $n_1=n$ in this case.

Suppose the contrary.
There exists $n'\ge n$ such that
$Q_{6+n'}$ is a node of $\Gamma_{6+n'}$ but
$Q_{6+i}$ is not for $i=n,\ldots,n'-1$.
If $n'=n$, then $E_{5+n}^2\le-3$ on $V_r$.
We infer that $h\circ\sigma_h$ cannot contract $E_{5+n}$ in this case.
Thus $n'>n$.
The weighted dual graph of $\Gamma_{7+n'}$ has the following shape.
\begin{center}
{\scalebox{1}{\includegraphics{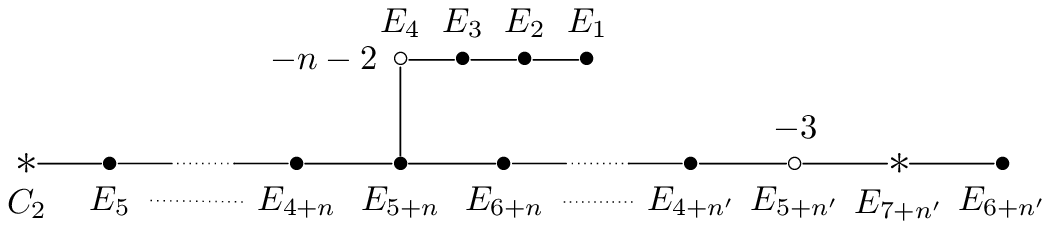}}}
\end{center}
The remaining blowings-up do not change curves other than
$E_{5+n'}$, $E_{6+n'}$ and $E_{7+n'}$.
It follows that
$h\circ\sigma_h$ first contracts $C_2,E_5,\ldots,E_{4+n'}$
in this order.
After that, $E_4$ must be a ($-1$)-curve.
Thus $n'=2n+1$.
If $r=7+n'$, then
$E_3^2=E_{5+n'}^2=-1$ after the contraction of $E_4$,
which is absurd.
For the same reason,
we have $r\ge10+n'$ and $Q_{i+n'}\in E_{i+n'}\cap E_{5+n'}$
for $i=7,8,9$.
If $r=10+n'$, then $h\circ\sigma_h$ cannot contract $E_{9+n'}$.
If $Q_{10+n'}\in E_{10+n'}\cap E_{5+n'}$,
then $h\circ\sigma_h$ cannot contract $E_{5+n'}$.
Thus $Q_{10+n'}\in E_{10+n'}\setminus E_{5+n'}$.
The situation on $V_{10+n'}$ is similar to that on $V_5$.
The curve $E_{5+n'}$ plays the role of $C_2$,
$E_{10+n'}$ does that of $E_5$ and so on.
\end{proof}

Let $C_2$ and $f$ be as in Theorem~\ref{thm:nonams}.
By applying Lemma~\ref{lem:y2} to $C_2$ and $f$,
we determine the minimal embedded resolution
$\sigma:V'\rightarrow\SP^2$ of the cusp of $C$.
On $\SP^2$, $f^{-1}$ has a single base point $Q\in C_2$.
We have either
$\{Q\}=C_2\cap f(C)$ or $Q\in C_2\setminus f(C)$.
By the same arguments as in Section~\ref{sec:pfif},
one can prove that
the weighted dual graph of $\sigma^{-1}(C)$
near $C'$ has the following shape.
In particular, we have $(C')^2=3$.
\begin{center}
\scalebox{1}{\includegraphics{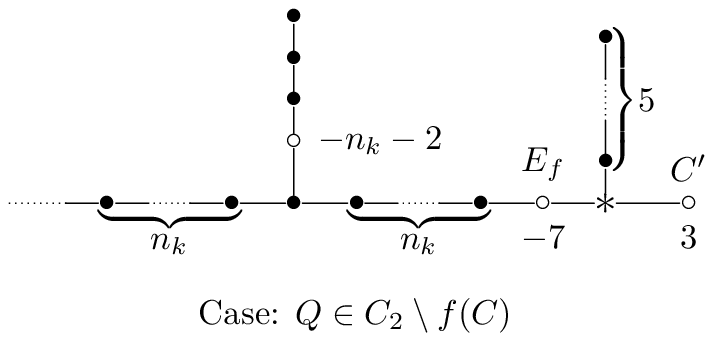}}
\end{center}
\begin{center}
\scalebox{1}{\includegraphics{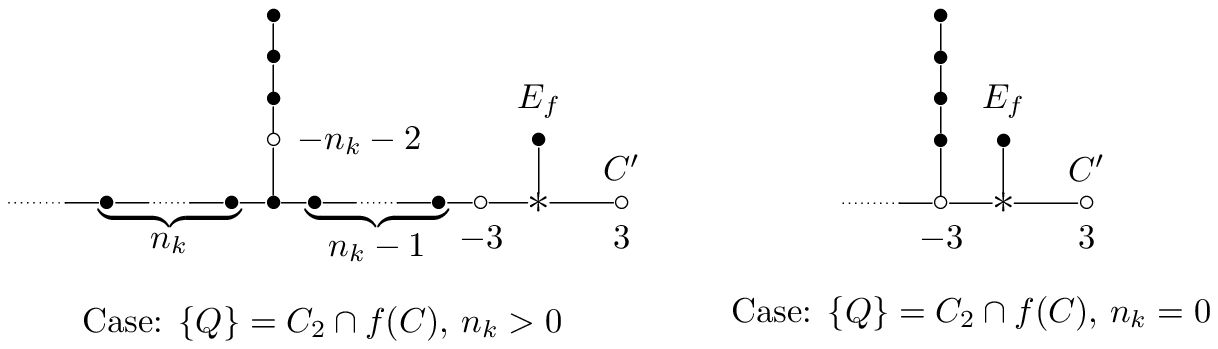}}
\end{center}
Finally, we show the existence of curves $C$ with $(C')^2=3$.
We start with
a line $L$,
a smooth conic $C_2$ and a smooth cubic $C_3$
meeting each other as in (\ref{eqn1}).
The next lemma shows the existence of such curves.
\begin{lemma}\label{lem:cubic}
Let $L,C_2,C_3$ be plane curves defined by the equations:
$L:x=0$, $C_2:f_2:=xz-y^2=0$, $C_3:(ax+2by-z)f_2+x^3=0$, $(a+b^2)^2\ne-4$.
Then they are smooth,
$C_3\cdot L=R_1+2R_2$
and
$C_3\cdot C_2=6R_2$,
where $R_1=(0,1,2b)$, $R_2=(0,0,1)$.
\end{lemma}
\begin{proof}
One can show that $C_3$ is projectively equivalent to
the curve defined by $y^2z=x^3+(a+b^2)x^2z-xz^2$.
It follows that
$C_3$ is smooth if and only if the right-hand side does not have
a multiple component.
\end{proof}

For given $k$, $n_1,\ldots,n_k$ and $Q\in C_2$,
let $h$ be a birational map given in Lemma~\ref{lem:y2}
with  the base point $Q$.
Let $C$ denote the strict transform of $C_3$
via $h$.
Then $C$ is an elliptic unicuspidal plane curve
satisfying the condition $(C')^2=3$.
\begin{example*}
Let $h$ be the birational map given in Lemma~\ref{lem:y2}
with $k=1$, $n_k=0$.
We may assume that
the base point of $h$ is $(0,0,1)$
and $C_2$ is defined by $f_2:=xz-y^2=0$.
Then
$h=h^{-1}=(xf_2^2,-f_3f_2,f_5)$,
where
$f_3=(cx+y)f_2+x^3$,
$f_5=(2x^2(cx+y)+(c^2x+2cy+z)f_2)f_2+x^5$.
Let $C_3$ be a smooth cubic meeting with $C_2$ at a single point $Q$
and
$C$ the strict transform of $C_3$ via $h$.
\begin{enumerate}
\item[(i)] $Q=(0,0,1)$.
The curve $C$ is a quintic having five double points at $(0,0,1)$
and is defined by
$axf_2^2-2bf_3f_2-f_5+x^3f_2=0$, $(a+b^2)^2\ne-4$.
\item[(ii)] $Q\ne(0,0,1)$.
The curve $C$ is of degree fifteen
and has six singular points at $(0,0,1)$
of multiplicity six.
It is defined by
$(af_5-2bf_3f_2-xf_2^2)f_2^5+f_5^3=0$, $(a+b^2)^2\ne-4$.
\end{enumerate}
\end{example*}
\begin{proof}
Let $L_x$, $L_y$ and $L_z$ denote the lines defined by
$x=0$, $y=0$ and $z=0$, respectively.
Put $g=h^{-1}$.
We have
$g^{\ast}(L_x)=L_x+2C_2$,
$g^{\ast}(L_y)=N+C_2$ and
$g^{\ast}(L_z)=C_5$.
Here $N$ (resp.~$C_5$) is a nodal cubic (unicuspidal quintic)
such that
$C_5$ has six double points at $P_1=(0,0,1)$
and that
$N C_2=6P_1$, $C_5 L_x=4 P_1+P_2$, $C_5 C_2=10 P_1$, $C_5 N=15 P_1$,
where $P_2\in L_x\setminus\{P_1\}$.
One can show that $N$ is defined by
$f_3=0$ and has the parameterization
$\SP^1\ni(s,t)\mapsto(st^2,s^2t(s-ct),s(s-ct)^2-t^3)\in N$.
By using the latter, one can deduce that $C_5$ is defined by $f_5=0$.
Since $zf_2^{12}|f_5\circ g$,
one has $g=(xf_2^2,-f_3f_2,f_5)=h$.
If $Q=(0,0,1)$, then $C$ is defined by
$f\circ g/f_2^5=0$,
where $f=(ax+2by-z)f_2+x^3$, $(a+b^2)^2\ne-4$.
Otherwise we may assume $Q=(1,0,0)$.
The curve $C$ is defined by
$f(f_5,-f_3f_2,xf_2^2)=0$.
\end{proof}
\begin{acknowledgment}
The author would like to express his thanks to Professor Fumio Sakai
for his helpful comments  on a preliminary version of this paper.
\end{acknowledgment}

\ \\
\textsc{%
{\small
Department of Mathematics,
Faculty of Science,
Saitama University,\\
Shimo-Okubo 255,
Urawa Saitama 338--8570,
Japan.}}\\
{\small
\textit{E-mail address}:  \texttt{ktono@rimath.saitama-u.ac.jp}
}
\end{document}